\documentclass{article}
\usepackage[cp1251]{inputenc}
\usepackage[russian,english]{babel}
\usepackage{amsmath}
\usepackage{amssymb}
\usepackage{amsfonts}
\usepackage{amsthm}
\usepackage{url}

\def\thtext#1{
  \catcode`@=11
  \gdef\@thmcountersep{. #1}
  \catcode`@=12
}

\def\threst{
  \catcode`@=11
  \gdef\@thmcountersep{.}
  \catcode`@=12
}

 \catcode`@=11
 \def\.{.\spacefactor\@m}
 \catcode`@=12

\theoremstyle{definition}
\newtheorem{rk}{Remark}

\newcommand{\cF}{\mathcal{F}}

\renewcommand{\a}{\alpha}
\newcommand{\dl}{\delta}

\newcommand{\G}{\Gamma}

\newcommand{\s}{\sigma}

\newcommand{\diam}{\operatorname{diam}}

\newcommand{\0}{\emptyset}
\renewcommand{\:}{\colon}

\renewcommand{\ss}{\subset}
\newcommand{\x}{\times}
\def\rom#1{\emph{#1}}
\def\({\rom(}
\def\){\rom)}

\begin{document}
\title{Hausdorff Measure: Lost in Translation}
\author{Alexey A. Tuzhilin}
\date{}
\maketitle

\begin{abstract}
In the present article we describe how one can define Hausdorff measure allowing empty elements in coverings, and using infinite countable coverings only. In addition, we discuss how the use of different nonequivalent interpretations of the notion ``countable set'', that is typical for classical and modern mathematics, may lead to contradictions.
\end{abstract}

The $m$-dimensional Hausdorff measure is one of the main fundamental notions in Geometric Measure Theory. To define it, one needs to make an agreement what is the value of $0^0$. Also, there is a temptation to use the empty set as element of coverings, like it is naturally done in general measure theory. However, if we do that, we get the following problem in the case $m=0$: the diameters of the empty set and a singleton are equal to each other (and equal to $0$), but they have to contribute differently to the definition of zero-dimensional Hausdorff measure. Namely, the empty set should not influence on the value of the measure, but each singleton has to add $1$. This discrepancy can not be overcome by the agreement on $0^0$ only. The most popular approach follows the famous monograph by Federer ``Geometric Measure Theory''~\cite{Federer} in which the empty sets are forbidden to be elements of coverings, and the coverings can be not only infinite countable, but finite as well. In particular, the empty covering is allowed too. In the present article we describe how one can define Hausdorff measure allowing empty elements in coverings, and using infinite countable coverings only. In addition, we discuss how the use of different nonequivalent interpretations of the notion ``countable set'', that is typical for classical and modern mathematics, may lead to contradictions.

To start with, we recall the definition of Hausdorff measure, see~\cite{Federer}, pp. 169--171. We present slightly modified its version to reach the main subject faster by omitting the details unnecessary for us.

\section{Definition of Hausdorff measure}\label{sec:1}
\markright{\thesection.~Definition of out Hausdorff measure}
Let $X$ be an arbitrary metric space. Denote by $|xy|$ the distance between $x,y\in X$. Let $A$ be any \textbf{nonempty} subset of $X$. Then the \emph{diameter of $A$} is the value
\begin{equation}\label{eq:diam}
\diam A=\sup\bigl\{|xy|:x,y\in A\bigr\}.
\end{equation}

Choose
\begin{enumerate}
\item an arbitrary family $\cF$ of subsets of $X$, and
\item any function $\zeta\:\cF\to[0,\infty]$.
\end{enumerate}
For each $A\ss X$ and $\dl\in(0,\infty]$, a subfamily $G\ss\cF$ is called a \emph{$\dl$-covering of $A$} if $A\ss\cup_{S\in G}S$ and $\diam S<\dl$ for all $S\in G$.\footnote{In~\cite{Federer} the condition $\diam S<\dl$ is changed to $\diam S\le\dl$. This is not essential for Hausdorff measure definition because this measure does not depend on the sets of infinite diameter. However, if we decide to keep the condition from~\cite{Federer}, we will need to set $\infty^0$.} We put
\begin{equation}\label{eq:2}
\phi_\dl(A)=\inf\Bigl\{\sum_{S\in G}\zeta(S):\text{$G$ is a \textbf{countable} $\dl$-covering $A$}\Bigr\}. \\
\end{equation}
Since for $0<\dl<\s\le\infty$ we have $\phi_\dl\ge\phi_\s$, then for each $A\ss X$ we have
\begin{equation}\label{eq:3}
\psi(A):=\sup_{\dl>0}\phi_\dl(A)=\lim_{\dl\to0+}\phi_\dl(A).
\end{equation}

The way to obtain the function $\psi\:2^X\to[0,\infty]$ which we described above is called \emph{Caratheodory Construction}.

For any $m\in[0,\infty)$, the \emph{$m$-dimensional Hausdorff measure\/} can be obtained as a particular case of the Caratheodory Construction: as $\cF$ one takes the family of all \textbf{nonempty\/} subsets $S\ss X$; the function $\zeta$ is defined as
$$
\zeta(S)=\frac{\a(m)}{2^m}(\diam S)^m,\ \ \a(m)=\frac{\G(1/2)^m}{\G(1+m/2)},
$$
where $\G$ is the standard gamma-function:
$$
\G(t)=\int_0^\infty x^{t-1}e^{-x}\,dx.
$$

In the case of $m$-dimensional Hausdorff measure, the functions $\phi_\dl$ and $\psi$ are usually denoted as $H^m_\dl$ and $H^m$, respectively.

Then~\cite{Federer} states that $H^0$ is the counting measure (i.e., it corresponds to each finite subset its cardinality, and to each infinite one the symbol $\infty$). Notice that~\cite{Federer} does not discuss what is the value $(\diam S)^0$ for a singleton $S$.

\section{What is countable?}
\markright{\thesection.~What is countable?}
To start with, let us mention the difference in definitions of a countable set in Russian and English Wikipedia.

Russia Wikipedia~\cite{WikiRus} declares:
\begin{quote}
In the set theory, \emph{a countable set\/} is an infinite set whose elements can be enumerated by natural numbers.
\end{quote}
However, in the next paragraph it is written: ``Sometimes countable sets may be defined as the sets of the same cardinality as subsets of the set of natural numbers, in particular, countable sets can be finite''.

Now, let us turn to English Wikipedia~\cite{WikiEng}:
\begin{quote}
In mathematics, a \emph{countable set\/} is a set with the same cardinality (number of elements) as some subset of the set of natural numbers.
\end{quote}
Here we also see an addendum: ``Some authors use countable set to mean countably infinite alone. To avoid this ambiguity, the term at most countable may be used when finite sets are included and countably infinite, enumerable, or denumerable otherwise.''

Now, let us look in some popular textbooks and monographs. One of very popular in Russia classical textbook by Kolmogorov and Fomin~\cite{Kolmog}, page~23, states:
\begin{quote}
By \emph{a countable set\/} we call any set whose elements can be bijectively correspond to all natural numbers.
\end{quote}

Textbook of Rudin~\cite{Rudin}, page~25:
\begin{quote}
$A$ is \emph{countable\/} if $A\sim J$.
\end{quote}
Here $\sim$ denotes that there is a bijective correspondence, and $J$ the set of all positive integers.

Hausdorff~\cite{Hausdorff}, page~29:
\begin{quote}
Sets of this [i.e., $\aleph_0$] cardinality ... are all countable (or enumerable or denumerable).
\end{quote}
In a footnote in the same page it is written: ``Sets that are finite (including the null set) or countable will be called at most countable''.

Kuratowski and Mostowski~\cite{Kuratowski}, page~174:
\begin{quote}
A set is said to be countable (or denumerable) if it is either finite or equipollent with the set of all natural numbers.
\end{quote}

\textbf{In what follows, to distinguish strictly between the two meanings of the term ``countable'', the sets for which there is a bijection with the set of all natural numbers we call \emph{infinite countable}, and the sets for which there is a bijection with a subset of the set of all natural numbers we call \emph{at most countable}}.

\section{Misunderstandings}
\markright{\thesection.~Misunderstandings}
The ambiguity in the definition of ``countable set'' described above leads to misunderstandings. For example, if we assume that the notion ``countable'' means infinite countable, then the definition of Hausdorff measure given in~\cite{Federer} creates contradictions. Let us describe some of them.

Let $X$ be a nonempty metric space consisting of $n$ points. We are interesting in the counting measure $H^0$. Choose $\dl>0$ to be less than the least distance between different points of $X$. Then each $\dl$-covering $G$ is the family of singletons. Each singleton $S\in G$ contributes in the definition of $H^0_\dl$ by $(\diam S)^0$. On the other hand, $G$ has to be infinite countable by definition, thus, $H^0_\dl$ equals the sum of infinitely many terms of the same nonnegative value. Therefore, $H^0(X)$ may be equal either $0$, or $\infty$, thus $H^0$ is not the counting measure.

In account, we come to conclusion: \emph{the claim that $G$ is infinite countable, together with the claim that $\cF$ consists of nonempty sets, lead to contradiction}.

Now, let us try to withdraw the nonemptiness condition on the elements of $\cF$, but preserve the infinite countability condition on the families $G$. In this case we need to define what is the diameter of an empty set. There is a well-know approach to do that.

Namely, extend the standard definition of $\inf R$ and $\sup R$ given for nonempty $R\ss[0,\infty]$ to the case $R=\0$ as follows:
$$
\inf\0=\infty\ \text{and}\ \sup\0=0.
$$
Then the formula~(\ref{eq:diam}) becomes meaningful for empty sets $A$ too: $\diam\0=0$.

Let $A$ be the empty subset of a nonempty finite metric space $X$. Choose $\dl>0$ to be less than the least distance between different points of $X$. Then each covering $G$ from the definition of $H^0_\dl(A)$ consists of singletons $\{x\}\ss X$ and empty sets $\0\ss X$. Since $(\diam\0)^0=\bigl(\diam\{x\}\bigr)^0=0^0$, then we again obtain the sum of infinite number of the same nonnegative values, a contradiction.

In account, we have shown that the \emph{withdraw from nonemptiness condition of the $\cF$'s elements while the coverings $G$ remain infinitely countable, leads to contradiction again}.

Thus, we see that various attempts to interpret the notion of countable set as infinite countable lead to contradictions. To save the situation, we will use the term ``countable'' in the sense ``at most countable''. Then the definition of Hausdorff measure from~\cite{Federer} turns out to be correct. Recall that Federer forbids empty sets in coverings $G$, allows to these coverings to be at most countable, in particular, to be empty, and sets (implicitly) $0^0=1$ to use this agreement in definition of zero-dimensional Hausdorff measure.

\section{Not only we}
\markright{\thesection.~Not only we}
It seems that the problems described above bother not only us. For instance, Halmos in~\cite{Halmos}, p. 53, defines $H^m$ by the formulas~(\ref{eq:2}) and~(\ref{eq:3}) for $m>0$ only. With this approach, we can set $H^0$ to be the counting measure just \textbf{by definition}.

Another attempt to overcome various contradictions (under ``wrong'' understanding the term countable) was undertaken in~\cite{Eilenberg} by Eilenberg and Harrold. They wrote: \emph{it is agreed that $(\diam A)^0=0,1$ according as $A$ is vacuous or not}. Of course, if one considers the $\diam A$ as a value from the segment $[0,\infty]$, then this definition looks not very rigorous. Anyway, it contains an explicit allusion to how one can resolve the contradiction in ``accurate'' way by preserving the solely infinite coverings, and allowing empty sets in them.

In the next section we present our definition of Hausdorff measure which was induced by~\cite{Eilenberg}. In our opinion, the necessity to consider infinite countable and finite coverings separately makes often the proofs unreasonably cumbersome, that is why our approach might be useful.

\section{Our definition of Hausdorff measure}
\markright{\thesection.~Our definition of Hausdorff measure}
To start with, we need to change slightly the Caratheodory Construction. Namely, we discard the condition that each $S\in\cF$ is nonempty and, instead of that, we put $\cF=2^X$. Also we replace $\zeta\:\cF\to[0,\infty]$ with the function $\zeta\:2^X\x[0,\infty)\to[0,\infty]$ possessing the following properties:
\begin{enumerate}
\item $\zeta(\0,m)=0$ for any $m\in[0,\infty)$;
\item $\zeta\bigl(\{x\},0\bigr)=1$ for any $x\in X$;
\item $\zeta\bigl(\{x\},m\bigr)=0$ for any $x\in X$ and any $m\in(0,\infty)$;
\item $\zeta\bigl(S,m\bigr)=(\a_m/2^m)(\diam S)^m$ for any $S\ss X$ consisting of at least $2$ elements, and any $m\in[0,\infty)$.
\end{enumerate}

For any $A\ss X$ and $\dl\in(0,\infty]$, a family $G$ of arbitrary (perhaps, empty) subsets of the space $X$ is called a \emph{$\dl$-covering of $A$} if $A\ss\cup_{S\in G}S$ and $\diam S<\dl$ for all $S\in G$.

Now, we can apply the modified Caratheodory Construction to get the desired result. Namely, we put
\begin{align}
&\label{align:Hk:1}H^m_\dl(A)=\inf\Bigl\{\sum_{S\in G}\zeta(S,m):\text{$G$ is an \textbf{infinite countable\/} $\dl$-covering of $A$}\Bigr\},\\
&\label{align:Hk:2}H^m(A)=\sup_{\dl>0}H^m_\dl(A)=\lim_{\dl\to0+}H^m_\dl(A).
\end{align}

In account, our definition preserves the both infinite countability of coverings $G$ and that $H^0$ is the counting measure. The empty subset can be covered by infinite family of empty subsets, and the condition $\zeta(\0,m)=0$ for all $m\in[0,\infty)$ gives us $H^m(\0)=0$.

\begin{rk}\label{rk:finiteReduction}
Since $\zeta(S,m)=0$ for $S=\0$, the sum in the right-hand side of formula~(\ref{align:Hk:1}) remains the same if we remove such $S$. This observation enables us to define $H^m_\dl$ by means of the both infinite and finite $\dl$-coverings $G$. Also, taking into account that for any $f\:G\to[0,\infty]$ we have $\sum_Gf=0$ for empty $G$, we can consider empty $\dl$-coverings $G$ as well (for empty $A$).

Thus, we obtain the following equivalent definition of Hausdorff measure: for any $\dl\in(0,\infty]$, $m\in[0,\infty)$, and $A\ss X$ we put
\begin{align*}
&H^m_\dl(A)=\inf\Bigl\{\sum_{S\in G}\zeta(S,m):\text{$G$ is \textbf{at most countable\/} $\dl$-covering of $A$}\Bigr\},\\
&H^m(A)=\sup_{\dl>0}H^m_\dl(A)=\lim_{\dl\to0+}H^m_\dl(A).
\end{align*}
\end{rk}

\begin{rk}
In our opinion, the function $\zeta(A,m)$ looks slightly artificial. It would be nice to find a natural interpretation of the function.
\end{rk}

\begin{rk}
It is easy to see that for $m>0$ our definition is equivalent to
\begin{align*}
&H^m_\dl(A)=\inf\Bigl\{\frac{\a_m}{2^m}\sum_{S\in G}(\diam S)^m:\text{$G$ is an infinite countable $\dl$-covering of $A$}\Bigr\},\\
&H^m(A)=\lim_{\dl\to0+}H^m_\dl(A)=\sup_{\dl>0}H^m_\dl(A),
\end{align*}
or, according to Remark~\ref{rk:finiteReduction},
\begin{align*}
&H^m_\dl(A)=\inf\Bigl\{\frac{\a_m}{2^m}\sum_{S\in G}(\diam S)^m:\text{$G$ is at most countable $\dl$-covering of $A$}\Bigr\},\\
&H^m(A)=\lim_{\dl\to0+}H^m_\dl(A)=\sup_{\dl>0}H^m_\dl(A).
\end{align*}
Here we do not need the function $\zeta$, and the only agreement we have to use is $\diam\0=0$ to define $\dl$-coverings which can contain empty sets. If we are interesting in $H^m$ for $m>0$ only, this variant (in fact, coinciding with the Halmos' one, see above) seems the most convenient.
\end{rk}

\date{\textbf{Acknowledgement: }The author thanks Hong Van Le and  Giovanni Alberti for their interest to this work and fruitful discussions.}

\markright{References}

\end{document}